\documentclass[10pt,twocolumn,twoside,a4paper]{article}

\usepackage[hyphens]{url}
\usepackage{xurl}
\usepackage[hidelinks]{hyperref}

\usepackage[left=15mm, right=15mm, top=15mm, bottom=15mm]{geometry}
\usepackage{subfig}
\usepackage{flushend}
\usepackage{indentfirst}
\usepackage{graphics}
\usepackage{amsmath,amsthm,amssymb}
\usepackage{multirow}
\usepackage{graphicx}
\usepackage{epstopdf}
\usepackage{float}
\usepackage[numbers]{natbib}
\usepackage[font=footnotesize,labelfont=bf]{caption}
\usepackage[labelsep=period]{caption}
\usepackage{fancyhdr}
\usepackage{mathrsfs}
\usepackage{enumerate}

\graphicspath{{./figure/}}

\rmfamily

\fancypagestyle{plain}{
\fancyhead{}
\fancyhead[L]{}
\fancyhead[R]{}
}
\newtheorem{thm}{Theorem}[section]
\newtheorem{lem}[thm]{Lemma}
\newtheorem{pro}[thm]{Proposition}

\theoremstyle{definition}
\newtheorem{defn}[thm]{Definition}
\newtheorem{exa}[thm]{Example}
\theoremstyle{remark}
\newtheorem{rem}[thm]{Remark}

\newcommand{\N}{\mathbb{N}}
\newcommand{\Rp}{\mathbb{R}_+}
\newcommand{\Acal}{\mathcal{A}}
\newcommand{\Adcal}{\mathcal{A}_d}
\newcommand{\Phical}{\Phi}

\begin{document}

\title{\huge \textbf{Fixed Point Results via a Uniform Class of $\mathcal{A}_d$-Contractions}}

\twocolumn[
\begin{@twocolumnfalse}

\author{\textbf{Prasun Panthi}$^{1}$, \textbf{Dinesh Panthi}$^{2,*}$\\
\footnotesize $^{1}$Department of Mathematics and Computer Science, Wabash College, Crawfordsville, IN 47933, USA\\
\footnotesize $^{2}$Department of Mathematics, Valmeeki Campus, Nepal Sanskrit University, Kathmandu 44600, Nepal\\
\footnotesize $^{*}$Corresponding Author: dinesh.panthi@nsu.edu.np}

\date{}
\maketitle

\end{@twocolumnfalse}
]

\noindent \textbf{\large{Abstract}}\hspace{2pt}
This paper studies fixed point arguments for dislocated metric spaces, where the self-distance of a point is not required to be zero. The starting point is the class $\Acal$ of control functions introduced by Akram, Zafar and Siddiqui for $\Acal$-contractions. In ordinary presentations of this class, the contraction factor may be obtained only after a particular comparison of two nonnegative numbers has been fixed. For an iterative proof, however, the same factor must control all steps of the Picard orbit. We isolate this issue by introducing a uniform subclass $\Adcal$, in which the contractive implication is governed by a single constant $k_\alpha<1$ attached to the function $\alpha$. Within the sequential theory of dislocated metrics, we prove four fixed point results: a theorem for one self-map, a theorem for a countable family of self-maps, an integral-type version, and a common fixed point theorem for two maps associated with two compatible dislocated metrics. We also include a detailed separation example showing that $\Adcal$ is a proper subclass of $\Acal$ and that the omitted uniformity condition cannot be recovered from the original $\Acal$-axiom. The examples are chosen to distinguish the new uniform condition from standard contraction hypotheses and to demonstrate that the results apply both to elementary interval models and to a function-space model with a nonzero fixed point.\\

\noindent \textbf{\large{Keywords}}\hspace{2pt} Common fixed point, dislocated metric space, fixed point, integral type contraction, sequence of mappings, $\Adcal$-contraction\\

\noindent\hrulefill

\section{\Large{Introduction}}

Banach's contraction principle remains one of the basic mechanisms for proving existence and uniqueness of fixed points \cite{banach}. Much of the later literature develops this mechanism in two related ways: either the distance structure is weakened, or the contractive inequality is replaced by a more flexible control condition. Dislocated metric spaces belong to the first direction. They arise from the work of Hitzler and Seda on dislocated topologies and generalized metrics \cite{hitzler-seda,hitzler}; unlike ordinary metric spaces, they allow the self-distance $d(x,x)$ to be positive. This feature connects dislocated metrics with metric-like and partial-metric ideas initiated in domain-theoretic settings by Matthews and subsequently studied in fixed point theory \cite{matthews,amini,karapinar-salimi,karapinar-survey,pasicki,zeyada,shoaib}.

The second direction concerns contractive conditions. Akram, Zafar and Siddiqui introduced the class $\Acal$ of control functions in order to treat $\Acal$-contractions in a unified way \cite{akram}. Integral variants of such contractions were later developed by Saha and Dey \cite{saha-dey}, in a line of research that is also related to integral-type contractions of Branciari and Rhoades \cite{branciari,rhoades}. Common fixed point theorems for pairs and families of maps under generalized contractive hypotheses have also been studied extensively; see, for instance, \cite{jungck,hussain2012,nashine2011,panthi2016integral,panthi2016four,panthi2018meir}.

The point of departure in the present paper is a small but important logical issue in the original $\Acal$-framework. In the definition of $\Acal$, the assertion that $a\le kb$ for some $k\in[0,1)$ is a pointwise assertion. It does not ensure that the same number $k$ is available uniformly along an iterative sequence. In a Picard iteration, however, one needs a single decay factor in order to obtain a geometric estimate of the form
\[
        d(x_n,x_{n+1})\le k^n d(x_0,x_1).
\]
Without such a uniform $k$, the usual Cauchy argument may fail even when each individual step looks contractive. This observation motivates the subclass $\Adcal$, in which the control function carries a global contraction constant.

The paper uses only the sequential language of dislocated metric spaces. This choice is deliberate. In a dislocated metric space, convergence behaves differently from ordinary metric convergence: if $x_n\to x$, then necessarily $d(x,x)=0$. Thus the natural location of the eventual fixed point is the zero-self-distance core
\[
        Y=\{x\in X:d(x,x)=0\}.
\]
We record the elementary lemmas needed for this fact before proving the fixed point results.

\subsection*{Contribution}

The main contribution is the identification of the exact uniformity condition needed to make the Picard proof work for $\Acal$-type controls in the dislocated metric setting. More precisely, we make the following contributions.
\begin{enumerate}[(i)]
\item We introduce the class $\Adcal$ by strengthening the pointwise implication in the class $\Acal$ to a uniform implication with a constant depending only on the control function.
\item We prove that this strengthening is genuine by constructing an explicit function in $\Acal\setminus\Adcal$ whose Picard orbit has no uniform geometric decay.
\item We obtain fixed point theorems for a single self-map, a sequence of self-maps, integral-type contractions, and two self-maps governed by two compatible dislocated metrics.
\item We give examples showing both the separation from the older $\Acal$ condition and the applicability of the new condition in finite-dimensional and function-space settings.
\end{enumerate}

\subsection*{Organization of the paper}

Section~2 collects the sequential preliminaries. Section~3 recalls $\Acal$, introduces $\Adcal$, and proves the main theorems. Section~4 provides the examples, including the separation example. Section~5 discusses the role of the hypotheses, and Section~6 concludes the paper.

\section{\Large{Preliminaries}}

Throughout the paper, $\Rp=[0,\infty)$. We use the following standard terminology from the literature on dislocated metric spaces.

\begin{quote}
\small
\begin{defn}[Dislocated metric, {\cite{hitzler-seda,hitzler,pasicki}}]
Let $X$ be a non-empty set and let $d:X\times X\to\Rp$ be a function satisfying:
\begin{enumerate}[(i)]
\item $d(x,y)=d(y,x)$ for all $x,y\in X$;
\item $d(x,y)=0$ implies $x=y$;
\item $d(x,y)\le d(x,z)+d(z,y)$ for all $x,y,z\in X$.
\end{enumerate}
Then $d$ is called a \emph{dislocated metric} (or \emph{$d$-metric}) on $X$, and $(X,d)$ is called a \emph{dislocated metric space}.
\end{defn}
\end{quote}

\begin{quote}
\small
\begin{defn}[Cauchy sequence, {\cite{hitzler,pasicki,karapinar-survey}}]
A sequence $\{x_n\}$ in a dislocated metric space $(X,d)$ is called \emph{$d$-Cauchy} if for every $\varepsilon>0$ there exists $n_0\in\N$ such that
\[
d(x_m,x_n)<\varepsilon \qquad \text{for all }m,n\ge n_0.
\]
\end{defn}
\end{quote}

\begin{quote}
\small
\begin{defn}[Sequential convergence, {\cite{hitzler,pasicki}}]
A sequence $\{x_n\}$ in a dislocated metric space $(X,d)$ is said to \emph{converge with respect to $d$} to $x\in X$ if
\[
\lim_{n\to\infty} d(x_n,x)=0.
\]
\end{defn}
\end{quote}

\begin{quote}
\small
\begin{defn}[Completeness, {\cite{hitzler,pasicki,panthi2016four}}]
A dislocated metric space $(X,d)$ is called \emph{complete} if every $d$-Cauchy sequence converges in $d$.
\end{defn}
\end{quote}

\begin{quote}
\small
\begin{defn}[Sequential continuity, {\cite{hitzler,pasicki}}]
Let $(X,d)$ and $(Y,\rho)$ be dislocated metric spaces. A mapping $T:X\to Y$ is said to be \emph{sequentially continuous} at $x\in X$ if $x_n\to x$ in $(X,d)$ implies $Tx_n\to Tx$ in $(Y,\rho)$. If this happens at every point of $X$, then $T$ is called \emph{sequentially continuous}.
\end{defn}
\end{quote}

The next elementary facts are standard consequences of the sequential convention for dislocated metrics and are included to fix notation for the later proofs \cite{hitzler,pasicki,karapinar-survey}.

\begin{lem}
Limits in a dislocated metric space are unique.
\end{lem}

\begin{proof}
Suppose that $x_n\to x$ and $x_n\to y$. Then, by the triangle inequality,
\[
d(x,y)\le d(x,x_n)+d(x_n,y)
\]
for every $n\in\N$. Letting $n\to\infty$ yields $d(x,y)=0$, and hence $x=y$.
\end{proof}

\begin{lem}\label{lem:selfcauchy}
If $\{x_n\}$ is a $d$-Cauchy sequence in a dislocated metric space $(X,d)$, then
\[
\lim_{n\to\infty} d(x_n,x_n)=0.
\]
\end{lem}

\begin{proof}
Let $\varepsilon>0$. Since $\{x_n\}$ is $d$-Cauchy, there exists $n_0\in\N$ such that $d(x_m,x_n)<\varepsilon$ whenever $m,n\ge n_0$. Taking $m=n$, we obtain $d(x_n,x_n)<\varepsilon$ for every $n\ge n_0$. Thus $d(x_n,x_n)\to 0$.
\end{proof}

\begin{lem}\label{lem:limitzero}
If $x_n\to x$ in a dislocated metric space $(X,d)$, then $d(x,x)=0$.
\end{lem}

\begin{proof}
By the triangle inequality,
\[
d(x,x)\le d(x,x_n)+d(x_n,x)=2d(x_n,x)
\]
for every $n\in\N$. Letting $n\to\infty$ gives $d(x,x)=0$.
\end{proof}

\begin{rem}
Lemma~\ref{lem:limitzero} shows that a constant sequence $\{x,x,x,\dots\}$ converges to $x$ if and only if $d(x,x)=0$. The fixed points obtained below always have zero self-distance, so this feature does not interfere with the main arguments.
\end{rem}

\begin{pro}\label{pro:Ymetric}
Let $(X,d)$ be a dislocated metric space and let
\[
Y=\{x\in X:d(x,x)=0\}.
\]
If $Y\neq\varnothing$ and $D:Y\times Y\to\Rp$ is defined by $D(x,y)=d(x,y)$ for all $x,y\in Y$, then $D$ is a metric on $Y$.
\end{pro}

\begin{proof}
Symmetry and the triangle inequality are inherited from $d$. If $D(x,y)=0$, then $d(x,y)=0$, so $x=y$ by axiom~(ii) of a dislocated metric. Therefore $D$ is a metric on $Y$.
\end{proof}

\section{\Large{Main Results}}

\subsection{The Classes $\Acal$ and $\Adcal$}

We first recall the original class of control functions introduced in \cite{akram}, and then introduce the uniform class $\Adcal$.

\begin{quote}
\small
\begin{defn}[Class $\Acal$, {\cite{akram}}]
Let $\Acal$ be the class of all functions $\alpha:\Rp^3\to\Rp$ satisfying:
\begin{enumerate}[(i)]
\item $\alpha$ is continuous on $\Rp^3$;
\item whenever
\[
a\le \alpha(a,b,b)\;\text{ or }\;a\le \alpha(b,a,b)\;\text{ or }\;a\le \alpha(b,b,a),
\]
one has $a\le kb$ for some $k\in[0,1)$.
\end{enumerate}
A self-map $T$ on a metric space $(X,d)$ is called an \emph{$\Acal$-contraction} if there exists $\alpha\in\Acal$ such that
\[
d(Tx,Ty)\le \alpha\bigl(d(x,y),d(x,Tx),d(y,Ty)\bigr)
\]
for all $x,y\in X$.
\end{defn}
\end{quote}

\begin{defn}[The uniform class $\Adcal$]\label{def:Ad}
Let $\Adcal$ be the class of all functions $\alpha:\Rp^3\to\Rp$ satisfying:
\begin{enumerate}
\item[(A1)] $\alpha$ is continuous on $\Rp^3$;
\item[(A2)] there exists a constant $k_\alpha\in[0,1)$ such that whenever
\[
a\le \alpha(a,b,b)\;\text{ or }\;a\le \alpha(b,a,b)\;\text{ or }\;a\le \alpha(b,b,a),
\]
then $a\le k_\alpha b$.
\end{enumerate}
\end{defn}

\begin{rem}\label{rem:inclusion}
Every $\alpha\in\Adcal$ belongs to $\Acal$, since the uniform constant $k_\alpha$ in axiom~(A2) serves as the existential $k$ required by the original class. The inclusion $\Adcal\subset\Acal$ is strict: Example~\ref{exa:separation} exhibits an $\alpha\in\Acal\setminus\Adcal$. The extra requirement in Definition~\ref{def:Ad} is precisely the uniformity of the contractive factor, which is the feature needed in our iterative estimates so that the same geometric constant is valid at every step of the recurrence.
\end{rem}

\begin{pro}\label{pro:examplesalpha}
Let $0\le \lambda<1$. Then each of the following functions belongs to $\Adcal$:
\begin{enumerate}[(i)]
\item $\alpha_1(u,v,w)=\lambda u$;
\item $\alpha_2(u,v,w)=\dfrac{\lambda u}{1+v+w}$;
\item $\alpha_3(u,v,w)=\dfrac{\lambda u}{1+u+v+w}$.
\end{enumerate}
\end{pro}

\begin{proof}
Each function is continuous on $\Rp^3$.

For $\alpha_1$, if $a\le \alpha_1(a,b,b)=\lambda a$, then $a=0$, and hence $a\le \lambda b$. If $a\le \alpha_1(b,a,b)=\lambda b$ or $a\le \alpha_1(b,b,a)=\lambda b$, then again $a\le \lambda b$. Thus $\alpha_1\in\Adcal$ with $k_{\alpha_1}=\lambda$.

For $\alpha_2$, if $a\le \alpha_2(a,b,b)=\frac{\lambda a}{1+2b}$, then $a=0$, because $\lambda/(1+2b)<1$. Moreover,
\[
a\le \alpha_2(b,a,b)=\frac{\lambda b}{1+a+b}\le \lambda b,
\]
and similarly in the third case. Hence $\alpha_2\in\Adcal$ with $k_{\alpha_2}=\lambda$.

For $\alpha_3$, if $a\le \alpha_3(a,b,b)=\frac{\lambda a}{1+a+2b}$, then $a=0$. Also,
\[
a\le \alpha_3(b,a,b)=\frac{\lambda b}{1+a+2b}\le \lambda b,
\]
and the same estimate holds in the third case. Therefore $\alpha_3\in\Adcal$ with $k_{\alpha_3}=\lambda$.
\end{proof}

\begin{defn}[$\Adcal$-contraction]
A self-mapping $T$ on a dislocated metric space $(X,d)$ is called an \emph{$\Adcal$-contraction} if there exists $\alpha\in\Adcal$ such that
\begin{equation}\label{eq:adcontraction}
 d(Tx,Ty)\le \alpha\bigl(d(x,y),d(x,Tx),d(y,Ty)\bigr)
\end{equation}
for all $x,y\in X$.
\end{defn}

\subsection{Fixed Point Theorem for a Single Mapping}

The first theorem gives the basic Picard-iteration result for one self-map.

\begin{thm}\label{thm:single}
Let $(X,d)$ be a complete dislocated metric space and let $T:X\to X$ be an $\Adcal$-contraction. Then $T$ has a unique fixed point $z\in X$. Moreover, $d(z,z)=0$.
\end{thm}

\begin{proof}
Start with an arbitrary point $x_0\in X$ and generate its Picard orbit by setting $x_{n+1}=Tx_n$ for $n\ge0$. If for some index $n$ the distance $d(x_n,x_{n+1})$ is zero, axiom (ii) of the dislocated metric gives $x_n=x_{n+1}=Tx_n$, and a fixed point has already been found. We therefore continue under the nontrivial assumption that every successive distance is positive.

By \eqref{eq:adcontraction},
\begin{align*}
 d(x_n,x_{n+1})
 &=d(Tx_{n-1},Tx_n) \\
 &\le \alpha\bigl(d(x_{n-1},x_n),d(x_{n-1},x_n),d(x_n,x_{n+1})\bigr)
\end{align*}
for every $n\ge 1$. Since $\alpha\in\Adcal$, axiom~(A2) yields $k\in[0,1)$ such that
\[
d(x_n,x_{n+1})\le k\, d(x_{n-1},x_n)\qquad (n\ge 1).
\]
Iterating this one-step estimate gives
\begin{equation}\label{eq:geometric}
 d(x_n,x_{n+1})\le k^n d(x_0,x_1) \qquad (n\ge 0).
\end{equation}
Now take $m>n$. Chaining the points $x_n,x_{n+1},\ldots,x_m$ and then using \eqref{eq:geometric}, we obtain
\[
d(x_n,x_m)\le \sum_{j=n}^{m-1}d(x_j,x_{j+1})
\le \frac{k^n}{1-k}d(x_0,x_1)\to 0,
\]
so $\{x_n\}$ is $d$-Cauchy. Completeness of $(X,d)$ supplies a point $z\in X$ such that $x_n\to z$. Lemma~\ref{lem:limitzero} then places this limit in the zero-self-distance core, namely $d(z,z)=0$.

For each $n\in\N$,
\begin{align*}
 d(z,Tz)
 &\le d(z,x_{n+1})+d(Tx_n,Tz) \\
 &\le d(z,x_{n+1})+\alpha\bigl(d(x_n,z),d(x_n,x_{n+1}),d(z,Tz)\bigr).
\end{align*}
Letting $n\to\infty$ and using continuity of $\alpha$,
\[
d(z,Tz)\le \alpha(0,0,d(z,Tz)).
\]
Axiom~(A2) with $a=d(z,Tz)$ and $b=0$ gives $d(z,Tz)=0$, so $Tz=z$.

For the self-distance, $d(z,z)=d(Tz,Tz)\le\alpha(d(z,z),0,0)$, and (A2) forces $d(z,z)=0$ (already noted above).

It remains to exclude a second fixed point. Suppose $w=Tw$. Applying the preceding self-distance argument to $w$ gives $d(w,w)=0$. By \eqref{eq:adcontraction},
\[
d(w,z)=d(Tw,Tz)\le \alpha\bigl(d(w,z),0,0\bigr),
\]
and (A2) gives $d(w,z)=0$, hence $w=z$.
\end{proof}

\subsection{A Common Fixed Point Theorem for a Sequence of Mappings}

The next result treats a family of self-maps under one common control function. Results of this type are part of the common fixed point literature for families and pairs of mappings \cite{jungck,hussain2012,panthi2016four}.

\begin{thm}\label{thm:sequence}
Let $(X,d)$ be a complete dislocated metric space and let $\{f_n\}_{n=1}^{\infty}$ be a sequence of self-mappings on $X$. Suppose that there exists $\alpha\in\Adcal$ such that
\begin{equation}\label{eq:sequence}
 d(f_i x,f_j y)\le \alpha\bigl(d(x,y),d(x,f_i x),d(y,f_j y)\bigr)
\end{equation}
for all $x,y\in X$ and all $i,j\in\N$. Then $\{f_n\}_{n=1}^{\infty}$ has a unique common fixed point $z\in X$, and $d(z,z)=0$.
\end{thm}

\begin{proof}
Pick any $x_0\in X$ and build an orbit by applying the maps in order: $x_n=f_nx_{n-1}$ for $n\ge1$. By \eqref{eq:sequence} with $i=n$, $j=n+1$, $x=x_{n-1}$ and $y=x_n$,
\[
d(x_n,x_{n+1})\le \alpha\bigl(d(x_{n-1},x_n),d(x_{n-1},x_n),d(x_n,x_{n+1})\bigr).
\]
Axiom~(A2) yields $k\in[0,1)$ with $d(x_n,x_{n+1})\le k\,d(x_{n-1},x_n)$. The estimate is identical in form to the estimate used in Theorem~\ref{thm:single}; summing the resulting geometric tail shows that $\{x_n\}$ is $d$-Cauchy. Hence $x_n\to z$ for some $z\in X$, and Lemma~\ref{lem:limitzero} gives $d(z,z)=0$.

Fix $n\in\N$. For every $m\in\N$,
\begin{align*}
 d(z,f_n z)
 &\le d(z,x_{m+1})+d(f_{m+1}x_m,f_n z) \\
 &\le d(z,x_{m+1})\\
 &\quad+\alpha\bigl(d(x_m,z),d(x_m,x_{m+1}),d(z,f_n z)\bigr).
\end{align*}
Letting $m\to\infty$ and using continuity of $\alpha$,
\[
d(z,f_n z)\le \alpha(0,0,d(z,f_n z)).
\]
Axiom~(A2) gives $d(z,f_n z)=0$, so $z=f_n z$. Since $n$ is arbitrary, $z$ is a common fixed point, and $d(z,z)=0$ follows as in Theorem~\ref{thm:single}.

If $w$ is another common fixed point, then $d(w,w)=0$ and, by \eqref{eq:sequence},
\[
d(z,w)=d(f_i z,f_j w)\le \alpha\bigl(d(z,w),0,0\bigr)
\]
for any $i,j\in\N$. Axiom~(A2) gives $d(z,w)=0$, so $z=w$.
\end{proof}

\subsection{Integral Type Contractions}

Following the integral-type tradition in fixed point theory \cite{branciari,rhoades,saha-dey,panthi2016integral}, let $\Phical$ denote the following class.
\begin{quote}
\small
Let $\Phical$ be the class of all functions $\phi:[0,\infty)\to[0,\infty)$ which are Lebesgue integrable, summable on each compact subset of $[0,\infty)$, non-negative, and satisfy
\[
\int_0^{\varepsilon}\phi(t)\,dt>0\qquad \text{for every }\varepsilon>0.
\]
\end{quote}

\begin{pro}\label{pro:psi-metric}
Let $(X,d)$ be a dislocated metric space, and let $\Psi:[0,\infty)\to[0,\infty)$ be continuous, nondecreasing, subadditive, and such that $\Psi(t)=0$ if and only if $t=0$. Define $D_\Psi(x,y)=\Psi(d(x,y))$ for $x,y\in X$. Then $D_\Psi$ is a dislocated metric on $X$. Moreover, if $(X,d)$ is complete, then $(X,D_\Psi)$ is also complete.
\end{pro}

\begin{proof}
Symmetry is immediate. If $D_\Psi(x,y)=0$, then $\Psi(d(x,y))=0$, so $d(x,y)=0$ and $x=y$. Since $\Psi$ is nondecreasing and subadditive,
\[
D_\Psi(x,y)\le \Psi(d(x,z)+d(z,y))\le D_\Psi(x,z)+D_\Psi(z,y),
\]
so $D_\Psi$ is a dislocated metric.

Let $\{x_n\}$ be $D_\Psi$-Cauchy and let $\varepsilon>0$. Since $\Psi(\varepsilon)>0$, there exists $n_0$ with $D_\Psi(x_m,x_n)<\Psi(\varepsilon)$ for $m,n\ge n_0$. Monotonicity of $\Psi$ then forces $d(x_m,x_n)<\varepsilon$, so $\{x_n\}$ is $d$-Cauchy. By completeness there exists $x$ with $d(x_n,x)\to 0$, and continuity of $\Psi$ gives $D_\Psi(x_n,x)\to 0$.
\end{proof}

The following theorem adapts the preceding argument to an integral-type control condition.

\begin{thm}\label{thm:integral}
Let $(X,d)$ be a complete dislocated metric space and let $T:X\to X$ be a self-mapping. Suppose that there exist $\alpha\in\Adcal$ and $\phi\in\Phical$ such that $\Psi(s)=\int_0^s \phi(t)\,dt$ is subadditive on $[0,\infty)$ and
\begin{equation}\label{eq:integral}
\begin{aligned}
&\int_0^{d(Tx,Ty)}\phi(t)\,dt \;\le\; \alpha\!\left(\int_0^{d(x,y)}\phi(t)\,dt,\right.\\
&\qquad\left.\int_0^{d(x,Tx)}\phi(t)\,dt,\,\int_0^{d(y,Ty)}\phi(t)\,dt\right)
\end{aligned}
\end{equation}
for all $x,y\in X$. Then $T$ has a unique fixed point $z\in X$, and $d(z,z)=0$.
\end{thm}

\begin{proof}
Proposition~\ref{pro:psi-metric} allows us to replace $d$ by the composed dislocated metric $D_\Psi(x,y)=\Psi(d(x,y))$. This new space is complete. Condition~\eqref{eq:integral} reads $D_\Psi(Tx,Ty)\le\alpha(D_\Psi(x,y),D_\Psi(x,Tx),D_\Psi(y,Ty))$, so $T$ is an $\Adcal$-contraction on $(X,D_\Psi)$. Theorem~\ref{thm:single} yields a unique fixed point $z$, and $0=D_\Psi(z,z)=\Psi(d(z,z))$ implies $d(z,z)=0$.
\end{proof}

\subsection{A Common Fixed Point Theorem with Two Compatible Dislocated Metrics}

The final theorem uses two dislocated metrics on the same set. Such compatibility assumptions are standard in common fixed point arguments with more than one map or more than one distance structure \cite{jungck,panthi2016four,panthi2018meir}.

\begin{thm}\label{thm:two}
Let $(X,d)$ and $(X,\delta)$ be two dislocated metric spaces on the same set $X$. Let $S,T:X\to X$ be self-mappings such that:
\begin{enumerate}[(i)]
\item $d(x,y)\le \delta(x,y)$ for all $x,y\in X$;
\item $(X,d)$ is complete;
\item $T$ is sequentially continuous with respect to $d$;
\item there exists $\alpha\in\Adcal$ such that
\begin{align}
      \delta(Tx,Sy) &\le \alpha\bigl(\delta(x,y),\delta(x,Tx),\delta(y,Sy)\bigr), \label{eq:two1}\\
      \delta(Sy,Tx) &\le \alpha\bigl(\delta(y,x),\delta(y,Sy),\delta(x,Tx)\bigr) \label{eq:two2}
\end{align}
for all $x,y\in X$;
\item whenever $x_n\to x$ in $d$, one has $\delta(x_n,x)\to 0$.
\end{enumerate}
Then $S$ and $T$ have a unique common fixed point $z\in X$. Moreover, $d(z,z)=\delta(z,z)=0$.
\end{thm}

\begin{rem}
A convenient sufficient condition for hypothesis~(v) is the existence of $C>0$ such that $\delta(x,y)\le Cd(x,y)$ for all $x,y\in X$.
\end{rem}

\begin{proof}
Choose $x_0\in X$ and define $x_{2n+1}=Tx_{2n}$, $x_{2n+2}=Sx_{2n+1}$ for $n\ge 0$. By \eqref{eq:two1},
\begin{align*}
\delta(x_{2n+1},x_{2n+2})\le{} &\alpha\bigl(\delta(x_{2n},x_{2n+1}),\\
&\delta(x_{2n},x_{2n+1}),\delta(x_{2n+1},x_{2n+2})\bigr),
\end{align*}
and axiom~(A2) yields $k\in[0,1)$ with $\delta(x_{2n+1},x_{2n+2})\le k\,\delta(x_{2n},x_{2n+1})$. By \eqref{eq:two2} the same constant gives $\delta(x_{2n+2},x_{2n+3})\le k\,\delta(x_{2n+1},x_{2n+2})$. Hence $\delta(x_n,x_{n+1})\le k\,\delta(x_{n-1},x_n)$ for all $n\ge 1$. Summing the successive $\delta$-distances gives a convergent geometric tail; consequently $\{x_n\}$ is $\delta$-Cauchy. By~(i) it is also $d$-Cauchy, so there exists $z\in X$ with $x_n\to z$ in $d$. Hypothesis~(v) then gives $\delta(x_n,z)\to 0$, and Lemma~\ref{lem:limitzero} gives $d(z,z)=0$.

Since $x_{2n}\to z$ in $d$ and $T$ is sequentially continuous, $Tx_{2n}=x_{2n+1}\to Tz$ in $d$. But $x_{2n+1}\to z$ in $d$, so $Tz=z$. The triangle inequality gives $\delta(z,z)\le 2\delta(z,x_n)\to 0$, so $\delta(z,z)=0$. Using \eqref{eq:two1} with $x=y=z$,
\[
\delta(z,Sz)=\delta(Tz,Sz)\le \alpha(0,0,\delta(z,Sz)),
\]
and (A2) gives $\delta(z,Sz)=0$, so $Sz=z$.

For uniqueness, if $w=Tw=Sw$, then $\delta(w,w)=0$ by the same argument, and \eqref{eq:two1} gives $\delta(z,w)\le\alpha(\delta(z,w),0,0)$, so (A2) forces $\delta(z,w)=0$, hence $z=w$.
\end{proof}

\section{\Large{Examples}}

The examples below have two purposes. The first one separates the new class from the older class $\Acal$, while the remaining examples show how the hypotheses can be checked in concrete settings.

\begin{exa}[An $\Acal$-contraction that is not an $\Adcal$-contraction]\label{exa:separation}
Let $X=[0,\infty)$ and define
\[
d(x,y)=|x-y|+x+y \qquad (x,y\in X).
\]
The axioms of a dislocated metric are readily verified. Since $x,y\ge0$, the expression also simplifies to $d(x,y)=2\max\{x,y\}$. Define
\[
T(x)=\frac{x}{1+x}
\qquad \text{and} \qquad
\alpha(u,v,w)=\frac{2u}{2+u}.
\]

\noindent\textit{$\alpha\in\Acal$.} Continuity of $\alpha$ on $\Rp^3$ is immediate from the formula. If $a\le\alpha(a,b,b)=\frac{2a}{2+a}$, then $a=0$. If $a\le\alpha(b,a,b)=\frac{2b}{2+b}$, then $a\le\frac{2}{2+b}b$, and $\frac{2}{2+b}<1$ for every $b>0$. The same holds for $\alpha(b,b,a)$. Hence $\alpha\in\Acal$.

\noindent\textit{$\alpha\notin\Adcal$.} Assume, for contradiction, that the stronger uniform requirement defining $\Adcal$ were satisfied by some $k\in[0,1)$. Choosing $a=\frac{2b}{2+b}$ so that $a=\alpha(b,a,b)$, axiom~(A2) would require $\frac{2}{2+b}\le k$ for all $b>0$. But the left-hand side tends to $1$ as $b\to 0^+$, a contradiction.

\noindent\textit{$T$ satisfies the $\Acal$-type inequality.} Since $T$ is increasing and $d(x,y)=2\max\{x,y\}$,
\begin{align*}
d(Tx,Ty)&=\frac{2\max\{x,y\}}{1+\max\{x,y\}}=\frac{2d(x,y)}{2+d(x,y)}\\
&=\alpha\bigl(d(x,y),d(x,Tx),d(y,Ty)\bigr).
\end{align*}

\noindent\textit{No uniform geometric rate exists.} With $x_0=1$ and $x_{n+1}=Tx_n$, a direct induction gives $x_n=\frac{1}{n+1}$, so
\[
\frac{d(x_{n+1},x_{n+2})}{d(x_n,x_{n+1})}=\frac{n+1}{n+2}\to 1.
\]
Thus no single constant $k<1$ can satisfy $d(x_{n+1},x_{n+2})\le k\,d(x_n,x_{n+1})$ for all $n$. This computation is the promised obstruction: the map is controlled in the older pointwise sense, but its orbit cannot be controlled by one geometric ratio. That is exactly the gap closed by Definition~\ref{def:Ad}.
\end{exa}

\begin{exa}[Single mapping, sequence of mappings, and integral type]\label{exa:basic}
Let $X=[0,1]$ and $d(x,y)=|x-y|+x+y$. Then $(X,d)$ is a dislocated metric space. Let $\{x_n\}$ be $d$-Cauchy. By Lemma~\ref{lem:selfcauchy}, $2x_n=d(x_n,x_n)\to 0$, so $d(x_n,0)=2x_n\to 0$. Hence $(X,d)$ is complete with zero-self-distance core $Y=\{0\}$.

Let $\alpha(u,v,w)=\frac{1}{2}u\in\Adcal$ (Proposition~\ref{pro:examplesalpha}).

\textit{Theorem~\ref{thm:single}.} Define $T(x)=\frac{x}{4}$. Then
\begin{align*}
d(Tx,Ty)&=\tfrac{1}{4}d(x,y)\le\tfrac{1}{2}d(x,y)\\
&=\alpha\bigl(d(x,y),d(x,Tx),d(y,Ty)\bigr).
\end{align*}
The unique fixed point is $z=0$.

\textit{Theorem~\ref{thm:sequence}.} For each $n\in\N$ define $f_n(x)=\frac{x}{n+3}$. For $i,j\in\N$ and $x,y\in X$,
\begin{align*}
d(f_i x,f_j y)&\le\frac{2x}{i+3}+\frac{2y}{j+3}\le\frac{x+y}{2}\\
&\le\frac{1}{2}d(x,y)=\alpha\bigl(d(x,y),d(x,f_i x),d(y,f_j y)\bigr).
\end{align*}
The unique common fixed point is $z=0$.

\textit{Theorem~\ref{thm:integral}.} Take $\phi\equiv 1$, so $\Psi(s)=s$ (subadditive) and $\phi\in\Phical$. With $T(x)=\frac{x}{4}$, condition~\eqref{eq:integral} reduces to $d(Tx,Ty)=\frac{1}{4}d(x,y)\le\frac{1}{2}d(x,y)$. The unique fixed point is $z=0$.
\end{exa}

\begin{exa}[Two dislocated metrics]\label{exa:two}
Let $X=[0,1]$, $d(x,y)=|x-y|+x+y$, and $\delta(x,y)=2d(x,y)$. Both are dislocated metrics, $d\le\delta$, and $(X,d)$ is complete by Example~\ref{exa:basic}. Set $T(x)=S(x)=\frac{x}{4}$ and $\alpha(u,v,w)=\frac{1}{2}u\in\Adcal$. Then $T$ is sequentially continuous in $d$, and $\delta(x_n,x)=2d(x_n,x)\to 0$ whenever $x_n\to x$ in $d$. Moreover,
\begin{align*}
\delta(Tx,Sy)&=2d\!\left(\tfrac{x}{4},\tfrac{y}{4}\right)=\tfrac{1}{2}d(x,y)=\tfrac{1}{4}\delta(x,y)\\
&\le\tfrac{1}{2}\delta(x,y)=\alpha\bigl(\delta(x,y),\delta(x,Tx),\delta(y,Sy)\bigr),
\end{align*}
and the same bound gives \eqref{eq:two2}. Theorem~\ref{thm:two} yields the unique common fixed point $z=0$.
\end{exa}

\begin{exa}[A nonlinear $\Adcal$-function on a function space]\label{exa:functional}
Let $h(t)=t$ for $t\in[0,1]$ and
\[
X=\bigl\{f\in C([0,1]):\|f-h\|_\infty\le 1\bigr\},
\]
\[
d_h(f,g)=\|f-g\|_\infty+\|f-h\|_\infty+\|g-h\|_\infty.
\]
Then $d_h$ is a dislocated metric on $X$. If $\{f_n\}$ is $d_h$-Cauchy, then $2\|f_n-h\|_\infty=d_h(f_n,f_n)\to 0$ by Lemma~\ref{lem:selfcauchy}, so $f_n\to h$ uniformly and $d_h(f_n,h)\to 0$. Hence $(X,d_h)$ is complete.

Define
\[
T(f)=h+\tfrac{1}{16}(f-h),
\qquad
\alpha(u,v,w)=\frac{u}{2(1+v+w)}.
\]
By Proposition~\ref{pro:examplesalpha}, $\alpha\in\Adcal$ with $\lambda=\frac{1}{2}$. Since $T(f)-h=\frac{1}{16}(f-h)$,
\[
d_h(Tf,Tg)=\tfrac{1}{16}d_h(f,g).
\]
Also $d_h(f,Tf)=2\|f-h\|_\infty\le 2$ and $d_h(g,Tg)\le 2$, so
\begin{align*}
\alpha\bigl(d_h(f,g),d_h(f,Tf),d_h(g,Tg)\bigr)
&\ge\frac{d_h(f,g)}{2(1+2+2)}\\
&=\frac{1}{10}d_h(f,g).
\end{align*}
Since $\frac{1}{16}<\frac{1}{10}$, the $\Adcal$-contraction condition holds. Theorem~\ref{thm:single} applies, and solving $Tf=f$ gives the unique fixed point $f=h$, the nonzero function $h(t)=t$.
\end{exa}

\section{\Large{Discussion}}

The preceding results depend on three structural choices. First, the constant in Definition~\ref{def:Ad} is required to be attached to the control function itself, not to a particular step of an iteration. This is the feature that turns a one-step comparison into a geometric estimate. Second, the proofs are formulated in the sequential language of dislocated metrics. This avoids adding auxiliary topology while still respecting the fact that limits must lie in the zero-self-distance core. Third, the common fixed point theorem with two dislocated metrics uses an explicit compatibility condition: convergence in the complete metric $d$ must be strong enough to imply convergence with respect to $\delta$.

The examples also indicate the boundary of the theory. Example~\ref{exa:separation} shows that the old $\Acal$ condition can allow a sequence of successive ratios tending to one. Therefore, replacing $\Adcal$ by $\Acal$ in Theorem~\ref{thm:single} would not be justified by the proof given here. On the other hand, Proposition~\ref{pro:examplesalpha} shows that many familiar control functions already belong to $\Adcal$, so the extra uniformity condition does not make the class empty or artificial.

\section{\Large{Conclusions}}

We introduced a uniform class $\mathcal{A}_d$ of control functions obtained from the 
original class $\mathcal{A}$ by requiring a global contractive factor suitable for 
iterative arguments in dislocated metric spaces. The inclusion $\mathcal{A}_d \subsetneq 
\mathcal{A}$ is strict (Example~\ref{exa:separation}), and the strictness is directly 
tied to the failure of the uniform geometric estimate that underlies our proofs. Within 
the sequential framework adopted in this paper, we established fixed point theorems for 
a single self-map, a sequence of self-maps, an integral type contraction, and, under a 
compatibility assumption, a pair of maps acting on a set equipped with two dislocated metrics.

The examples clarify two complementary aspects of the theory. Example~\ref{exa:separation} 
shows that the original class $\mathcal{A}$ does not in general supply the uniform geometric 
estimate required by our proofs, thereby motivating the passage to $\mathcal{A}_d$. The 
remaining examples show that the theorems apply not only to elementary interval models 
but also to nonlinear controls and nontrivial fixed points in a function-space setting.

This gives a sequential framework for fixed point arguments in dislocated spaces in which the role of the global contraction factor is explicit. Possible extensions include asymmetric distance structures, such as dislocated quasi-metric spaces, and applications to existence and uniqueness problems for nonlinear integral control equations.

\section*{Acknowledgements}

The research of the first author was supported by a summer research grant from the Dill Fund at Wabash College, established through the generous support of Michael Dill '71. This support was instrumental in the completion of this work.

\noindent\hrulefill

\begingroup
\sloppy
\raggedright
\renewcommand\refname{REFERENCES}

\endgroup
\end{document}